\documentclass[12pt,a4paper]{article}

\usepackage[colour]{ajr}

\newcommand{\flops}{\,\mbox{flops}}

\begin{document}

\title{Fast and accurate multigrid solution of Poissons
equation using diagonally oriented grids}
\author{A.J.~Roberts\thanks{Dept. Maths \& Comput, University of
Southern Queensland, Toowoomba, Queensland~4352, \textsc{Australia}.
\protect\url{mailto:aroberts@usq.edu.au}}}
\maketitle

\begin{abstract}
We solve Poisson's equation using new multigrid algorithms that
converge rapidly.  The novel feature of the 2D and 3D algorithms are
the use of extra diagonal grids in the multigrid hierarchy for a much
richer and effective communication between the levels of the
multigrid.  Numerical experiments solving Poisson's equation in the
unit square and unit cube show simple versions of the proposed
algorithms are up to twice as fast as correspondingly simple
multigrid iterations on a conventional hierarchy of grids.
\end{abstract}

\tableofcontents

\section{Introduction}

Multigrid algorithms are effective in the solution of elliptic
problems and have found many applications, especially in fluid
mechanics \cite[e.g.]{Mavriplis97}, chemical reactions in flows
\cite[e.g.]{Sheffer98} and flows in porous media \cite{Moulton98}.
Typically, errors in a solution may decrease by a factor of 0.1 each
iteration \cite[e.g.]{Mavriplis97b}.  The simple algorithms I present
decrease errors by a factor of $0.05$ (see Tables~\ref{tbl:2d}
and~\ref{tbl:3dopt} on pages~\pageref{tbl:2d}
and~\pageref{tbl:3dopt}).  Further gains in the rate of convergence
may be made by further research.

Conventional multigrid algorithms use a hierarchy of grids whose grid
spacings are all proportional to $2^{-\ell}$ where $\ell$ is the level
of the grid \cite[e.g.]{Zhang98}.  The promising possibility I report
on here is the use of a richer hierarchy of grids with levels of the
grids oriented diagonally to other levels.  Specifically, in 2D I
introduce in Section~\ref{S2d} a hierarchy of grids with grid spacings
proportional to $2^{-\ell/2}$ and with grids aligned at $45^\circ$ to
adjacent levels, see Figure~\ref{Fgrid2d}
(p\pageref{Fgrid2d}).\footnote{This paper is best viewed and printed
in colour as the grid diagrams and the ensuing discussions are all
colour coded.} In 3D the geometry of the grids is much more
complicated.  In Section~\ref{S3d} we introduce and analyse a
hierarchy of 3D grids with grid spacings roughly $2^{-\ell/3}$ on the
different levels, see Figure~\ref{Famal} (p\pageref{Famal}).  Now
Laplace's operator is isotropic so that its discretisation is
straightforward on these diagonally oriented grids.  Thus in this
initial work I explore only the solution of Poisson's equation.
\begin{equation}
	\nabla^2 u=f\,.
	\label{Epois}
\end{equation}
Given an approximation $\tilde u$ to a solution, each complete
iteration of a multigrid scheme seek a correction $v$ so that
$u=\tilde u+v$ is a better approximation to a solution of Poisson's
equation~(\ref{Epois}).  Consequently the update $v$ has to
approximately satisfy a Poisson's equation itself, namely
\begin{equation}
	\nabla^2v=r\,,
	\quad\mbox{where}\quad
	r=f-\nabla^2\tilde u\,,
	\label{Evpois}
\end{equation}
is the residual of the current approximation.  The multigrid
algorithms aim to estimate the error $v$ as accurately as possible
from the residual $r$.  Accuracy in the ultimate solution $u$ is
determined by the accuracy of the spatial discretisation in the
computation of the residual $r$: here we investigate second-order and
fourth-order accurate discretisations \cite[e.g.]{Zhang98} but so far
only find remarkably rapid convergence for second-order
discretisations.

The diagonal grids employed here are perhaps an alternative to the
semi-coarsening hierarchy of multigrids used by Dendy
\cite{Dendy97} in more difficult problems.

In this initial research we only examine the simplest reasonable
V-cycle on the special hierarchy of grids and use only one Jacobi
iteration on each grid.  We find in Sections~\ref{SSsr2}
and~\ref{SSsr3} that the smoothing restriction step from one grid to
the coarser diagonally orientated grid is done quite simply.  Yet the
effective smoothing operator from one level to that a factor of 2
coarser, being the convolution of two or three intermediate steps, is
relatively sophisticated.  One saving in using these diagonally
orientated grids is that there is no need to do any interpolation.
Thus the transfer of information from a coarser to a finer grid only
involves the simple Jacobi iterations described in
Sections~\ref{SSjp2} and~\ref{SSjp3}.  Performance is enhanced within
this class of simple multigrid algorithms by a little over relaxation
in the Jacobi iteration as found in Sections~\ref{SSopt2}
and~\ref{SSopt3}.  The proposed multigrid algorithms are found to be
up to twice as fast as comparably simple conventional multigrid
algorithms.

\section{A diagonal multigrid for the 2D Poisson equation}
\label{S2d}

\begin{figure}[tbp]
	\centering
	\includegraphics{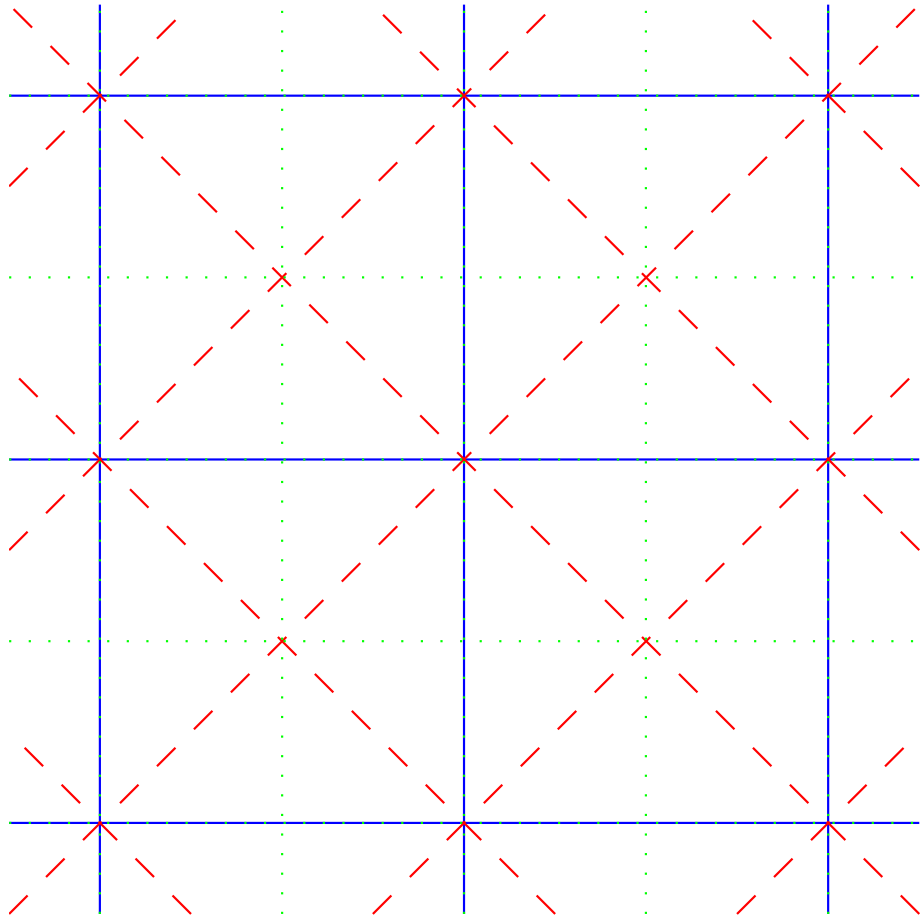} \caption{three levels of grids in the
	2D multigrid hierarchy: the dotted green grid is the finest,
	spacing $h$ say; the dashed red grid is the next finest diagonal
	grid with spacing $\sqrt2h$; the solid blue grid is the coarsest
	shown grid with spacing $2h$.  Coarser levels of the multigrid
	follow the same pattern.}
	\label{Fgrid2d}
\end{figure}
To approximately solve Poisson's equation~(\ref{Evpois}) in
two-dimensions we use a novel hierarchy of grids in the multigrid
method.  The length scales of the grid are $2^{-\ell/2}$.  If the
finest grid is aligned with the coordinate axes with grid spacing $h$
say, the first coarser grid is at $45^\circ$ with spacing $\sqrt2h$, the
second coarser is once again aligned with the axes and of spacing $2h$,
as shown in Figure~\ref{Fgrid2d}, and so on for all other levels on
the multigrid.  In going from one level to the next coarser level the
number of grid points halves.

\subsection{The smoothing restriction}
\label{SSsr2}

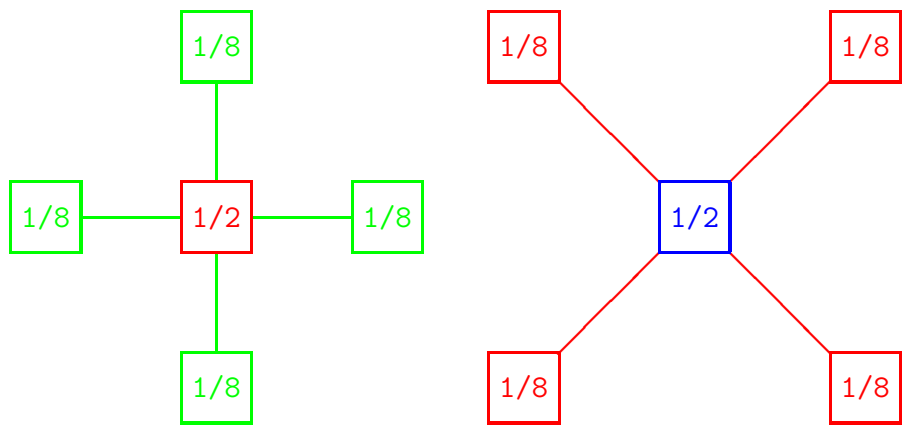
\begin{figure}[tbp]
	\centering
{\tt    \setlength{\unitlength}{0.25ex}
\begin{picture}(270,130)
\thicklines   {\color{red}
              \put(220,60){\line(1,-1){30}}
              \put(170,110){\line(1,-1){30}}
              \put(220,80){\line(1,1){30}}
              \put(170,30){\line(1,1){30}}
              \put(250,110){\framebox(20,20){1/8}}
              \put(250,10){\framebox(20,20){1/8}}
              \put(150,110){\framebox(20,20){1/8}}
              \put(150,10){\framebox(20,20){1/8}}
              }
              \put(200,60){\color{blue}\framebox(20,20){1/2}}
              {\color{green}
              \put(70,30){\line(0,1){30}}
              \put(30,70){\line(1,0){30}}
              \put(80,70){\line(1,0){30}}
              \put(70,80){\line(0,1){30}}
              \put(60,10){\framebox(20,20){1/8}}
              \put(110,60){\framebox(20,20){1/8}}
              \put(60,110){\framebox(20,20){1/8}}
              \put(10,60){\framebox(20,20){1/8}}
              }
              \put(60,60){\color{red}\framebox(20,20){1/2}}
\end{picture}}
\caption{restriction stencils are simple weighted averages of
	neighbouring grid points on all levels of the grid.}
	\label{Erest2}
\end{figure}
The restriction operator smoothing the residual from one grid to the next
coarser grid is the same at all levels.  It is simply a weighted
average of the grid point and the four nearest neighbours on the finer
grid as shown in Figure~\ref{Erest2}.  To restrict from a fine green
grid to the diagonal red grid
\begin{equation}
	r_{i,j}^{\ell-1}=\frac{1}{8}\left( 4r_{i,j}^\ell +r_{i-1,j}^\ell
	+r_{i,j-1}^\ell +r_{i+1,j}^\ell +r_{i,j+1}^\ell \right)\,,
	\label{Erest2r}
\end{equation}
whereas to restrict from a diagonal red grid to the coarser blue grid
\begin{equation}
	r_{i,j}^{\ell-1}=\frac{1}{8}\left( 4r_{i,j}^\ell +r_{i-1,j-1}^\ell
	+r_{i+1,j-1}^\ell +r_{i+1,j+1}^\ell +r_{i-1,j+1}^\ell \right)\,.
	\label{Erest2b}
\end{equation}
Each of these restrictions takes $6\flops$ per grid element.  Thus
assuming the finest grid is $n\times n$ with $N=n^2$ grid points, the
restriction to the next finer diagonal grid (red) takes approximately
$3N\flops$, the restriction to the next finer takes approximately
$3N/2\flops$, etc.  Thus to restrict the residuals up $\ell=2L$ levels
to the coarsest grid spacing of $H=2^Lh$ takes
\begin{equation}
	K_r\approx 6N\left(1-\frac{1}{4^L}\right)\flops \approx 6N\flops\,.
	\label{Ekrest2}
\end{equation}
In contrast a conventional nine point restriction operator from one
level to another takes $11\flops$ per grid point, which then totals to
approximately $3\frac{2}{3}N\flops$ over the whole conventional
multigrid hierarchy.  This is somewhat better than the proposed
scheme, but we make gains elsewhere.  In restricting from the green
grid to the blue grid, via the diagonal red grid, the restriction
operation is equivalent to a 17-point stencil with a much richer and
more effective smoothing than the conventional 9-point stencil.

\subsection{The Jacobi prolongation}
\label{SSjp2}

\begin{figure}[tbp]
	\centering
	\includegraphics[width=\textwidth]{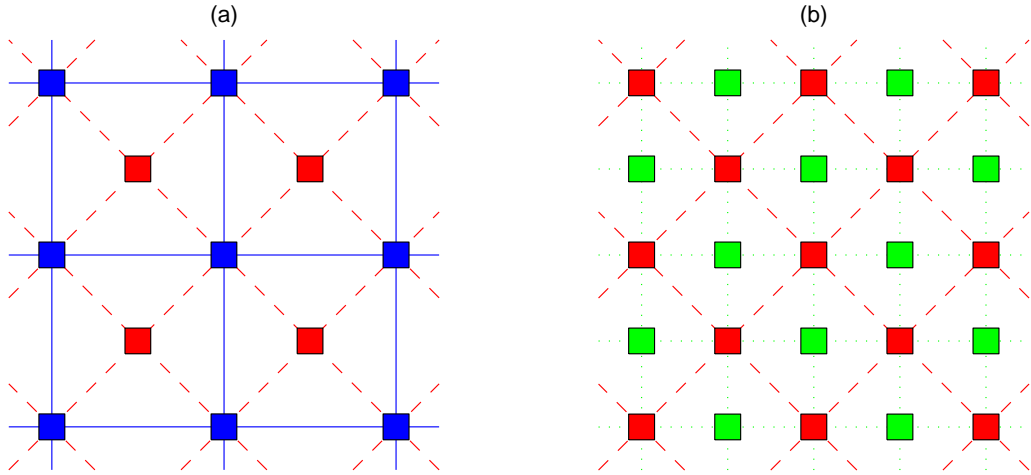}
	\caption{the interpolation in a prolongation step is replaced
	by simply a ``red-black'' Jacobi iteration: (a) compute the new
	values at the red grid points, then refine the values at the blue
	points; (b) compute the new values at the green points, then refine
	those at the red points.}
	\label{Fprol2d}
\end{figure}
One immediate saving is that there is no need to interpolate in the
prolongation step from one level to the next finer level.  For
example, to prolongate from the blue grid to the finer diagonal red grid,
shown in Figure~\ref{Fprol2d}(a), estimate the new value of $v$ at
the red grid points on level $\ell$ by the red-Jacobi iteration
\begin{equation}
	v_{i,j}^\ell=\frac{1}{4}\left( -2h^2r_{i,j}^\ell +v_{i-1,j-1}^{\ell-1}
	+v_{i+1,j-1}^{\ell-1} +v_{i+1,j+1}^{\ell-1} +v_{i-1,j+1}^{\ell-1}
\right)\,,
	\label{Ejacr}
\end{equation}
when the grid spacing on the red grid is $\sqrt2h$.  Then the values
at the blue grid points are refined by the blue-Jacobi iteration
\begin{equation}
	v_{i,j}^\ell=\frac{1}{4}\left( -2h^2r_{i,j}^\ell +v_{i-1,j-1}^\ell
	+v_{i+1,j-1}^\ell +v_{i+1,j+1}^\ell +v_{i-1,j+1}^\ell \right)\,.
	\label{Ejacb}
\end{equation}
A similar green-red Jacobi iteration will implicitly prolongate from
the red grid to the finer green grid shown in Figure~\ref{Fprol2d}(b).
These prolongation-iteration steps take $6\flops$ per grid point.
Thus to go from the red to the green grid takes $6N\flops$. As each
level of the grid has half as many grid points as the next finer, the
total operation count for the prolongation over the hierarchy from
grid spacing $H=2^Lh$ is
\begin{equation}
	K_p\approx 12N\left( 1-\frac{1}{4^L} \right)\flops \approx 12N\flops\,.
	\label{Ekprol2}
\end{equation}

The simplest (bilinear) conventional interpolation direct from the
blue grid to the green grid would take approximately $2N\flops$, to be
followed by $6N\flops$ for a Jacobi iteration on the fine green grid
(using simply $\nu_1=0$ and $\nu_2=1$).  Over the whole hierarchy this
takes approximately $10\frac{2}{3}N\flops$.  This is a little smaller
than that proposed here, but the proposed diagonal method achieves
virtually two Jacobi iterations instead of just one and so is more
effective.

\subsection{The V-cycle converges rapidly}

Numerical experiments show that although the operation count of the
proposed algorithm is a little higher than the simplest usual
multigrid scheme, the speed of convergence is much better.  The
algorithm performs remarkably well on test problems such as those in
Gupta et al \cite{Gupta97}.  I report a quantitative comparison
between the algorithms that show the diagonal scheme proposed here is
about twice as fast.

Both the diagonal and usual multigrid algorithms use $7N\flops$ to
compute the residuals on the finest grid.  Thus the proposed method
takes approximately $25N\flops$ per V-cycle of the multigrid
iteration, although 17\% more than the simplest conventional algorithm
that takes $21\frac{1}{3}N\flops$, the convergence is much faster.
Table~\ref{tbl:2d} shows the rate of convergence $\bar\rho_0\approx
0.1$ for this diagonal multigrid based algorithm.  The data is
determined using \matlab's sparse eigenvalue routine to find the
largest eigenvalue and hence the slowest decay on a $65\times 65$
grid.  This should be more accurate than limited analytical methods
such as a bi-grid analysis \cite{Ibraheem96}.  Compared with
correspondingly simple schemes based upon the usual hierarchy of
grids, the method proposed here takes much fewer iterations, even
though each iteration is a little more expensive, and so should be
about twice as fast.

\begin{table}[tbp]
	\centering
	\caption{comparison of cost, in flops, and performance for various
	algorithms for solving Poisson's equation in two spatial
	dimensions.  The column headed ``per iter'' shows the number of
	flops per iteration, whereas columns showing ``per dig'' are
	$\flops/\log_{10}\bar\rho$ and indicate the number of flops needed
	to compute each decimal digit of accuracy.  The right-hand columns
	show the performance for the optimal over relaxation parameter
	$p$.}
	\begin{tabular}{|l|r|lr|llr|}
		\hline
		algorithm & per iter & $\bar\rho_0$ & per dig & $p$ &
$\bar\rho$ &
		per dig  \\
		\hline
		diagonal, $\Ord{h^2}$ & $25.0N$ & .099 & $25.0N$ & 1.052 &
.052 & $19.5N$  \\
		usual, $\Ord{h^2}$ & $21.3N$ & .340 & $45.5N$ & 1.121 &
.260 & $36.4N$  \\
		\hline
		diagonal, $\Ord{h^4}$ & $30.0N$ & .333 & $62.8N$ & 1.200 &
.200 & $42.9N$  \\
		usual, $\Ord{h^4}$ & $26.3N$ & .343 & $56.6N$ & 1.216 &
.216 & $39.4N$  \\
		\hline
	\end{tabular}
	\label{tbl:2d}
\end{table}

Fourth-order accurate solvers in space may be obtained using the above
second-order accurate V-cycle as done by Iyengar \& Goyal
\cite{Iyengar90}.  The only necessary change is to compute the
residual $r$ in~(\ref{Evpois}) on the finest grid with a fourth-order
accurate scheme, such as the compact ``Mehrstellen'' scheme
\begin{eqnarray}
r_{i,j}&=&\frac{1}{12}\left( 8f_{i,j}
           +f_{i+1,j} +f_{i,j+1} +f_{i-1,j} +f_{i,j-1} \right)
           \nonumber\\&&{}
       -\frac{1}{6h^2}\left[ -20u_{i,j}
         +4\left(u_{i,j-1} +u_{i,j+1} +u_{i-1,j} +u_{i+1,j}\right)
           \right.\nonumber\\&&\quad\left.{}
	     +u_{i+1,j+1} +u_{i-1,j+1} +u_{i-1,j-1} +u_{i+1,j-1}
	      \right]\,.
	\label{Efos2}
\end{eqnarray}
Use the V-cycles described above to determine an approximate
correction $v$ to the field $u$ based upon these more accurate
residuals.  The operation count is solely increased by the increased
computation in the residual, from $7N\flops$ per iteration to
$12N\flops$ (the combination of $f$ appearing on the right-hand side
of~(\ref{Efos2}) need not be computed each iteration).  Numerical
experiments summarised in Table~\ref{tbl:2d} show that the multigrid
methods still converge, but the diagonal method has lost its
advantages.  Thus fourth order accurate solutions to Poisson's
equation are most quickly obtained by initially using the diagonal
multigrid method applied to the second order accurate computation of
residuals.  Then use a few multigrid iterations based upon the fourth
order residuals to refine the numerical solution.

\subsection{Optimise parameters of the V-cycle}
\label{SSopt2}

The multigrid iteration is improved by introducing a small amount of
over relaxation.

First we considered the multigrid method applied to the second-order
accurate residuals. Numerical optimisation over a range of
introduced parameter values suggested that the simplest, most robust
effective change was simply to introduce a parameter $p$ into the
Jacobi iterations~(\ref{Ejacr}--\ref{Ejacb}) to become
\begin{eqnarray}
	v_{i,j}^\ell&=&\frac{1}{4}\left( -2ph^2r_{i,j}^\ell
+v_{i-1,j-1}^{\ell-1}
	+v_{i+1,j-1}^{\ell-1} +v_{i+1,j+1}^{\ell-1} +v_{i-1,j+1}^{\ell-1}
\right)\,,
	\label{Eajacr}\\
	v_{i,j}^\ell&=&\frac{1}{4}\left( -2ph^2r_{i,j}^\ell +v_{i-1,j-1}^\ell
	+v_{i+1,j-1}^\ell +v_{i+1,j+1}^\ell +v_{i-1,j+1}^\ell \right)\,,
	\label{Eajacb}
\end{eqnarray}
on a diagonal red grid and similarly for a green grid.  An optimal
value of $p$ was determined to be $p=1.052$.  The parameter $p$ just
increases the weight of the residuals at each level by about 5\%.
This simple change, which does not increase the operation count,
improves the factor of convergence to $\bar\rho\approx 0.052$, which
decreases the necessary number of iterations to achieve a given
accuracy.  As Table~\ref{tbl:2d} shows, this diagonal multigrid is
still far better than the usual multigrid even with its optimal choice
for over relaxation.

Then we considered the multigrid method applied to the fourth-order
accurate residuals.  Numerical optimisation of the parameter $p$
in~(\ref{Eajacr}--\ref{Eajacb}) suggests that significantly more
relaxation is preferable, namely $p\approx 1.20$.  With this one
V-cycle of the multigrid method generally reduces the residuals by a
factor $\bar\rho\approx 0.200$.  This simple refinement reduces the
number of iterations required by about one-third in converging to the
fourth-order accurate solution.

\section{A diagonal multigrid for the 3D Poisson equation}
\label{S3d}

\begin{figure}[tbp]
	\centering
	\includegraphics[width=0.9\textwidth]{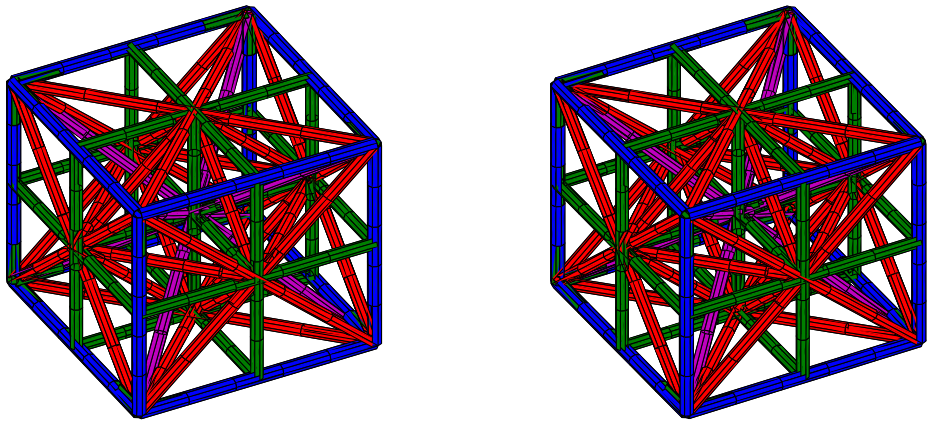}
\caption{one cell of an amalgam of four levels of the hierarchy of
grids used to form the multigrid V-cycle in 3D: green is the finest
grid shown; red is the next level coarser grid; magenta shows the next
coarser grid; and the blue cube is the coarsest to be shown.  This
stereoscopic view is to be viewed cross-eyed as this seems to be more
robust to changes of viewing scale.}
	\label{Famal}
\end{figure}

The hierarchy of grids we propose for solving Poisson's
equation~(\ref{Evpois}) in three-dimensions is significantly more
complicated than that in two-dimensions.  Figure~\ref{Famal} shows the
three steps between levels that will be taken to go from a fine
standard grid (green) of spacing $h$, via two intermediate grids (red
and magenta), to a coarser regular grid (blue) of spacing $2h$.  As we
shall discuss below, there is some unevenness in the hierarchy that
needs special treatment.

\subsection{The smoothing restriction steps}
\label{SSsr3}

The restriction operation in averaging the residuals from one grid to
the next coarser grid is reasonably straightforward.
\begin{itemize}
	\item
	\begin{figure}[tbp]
		\centering
		\includegraphics[width=0.9\textwidth]{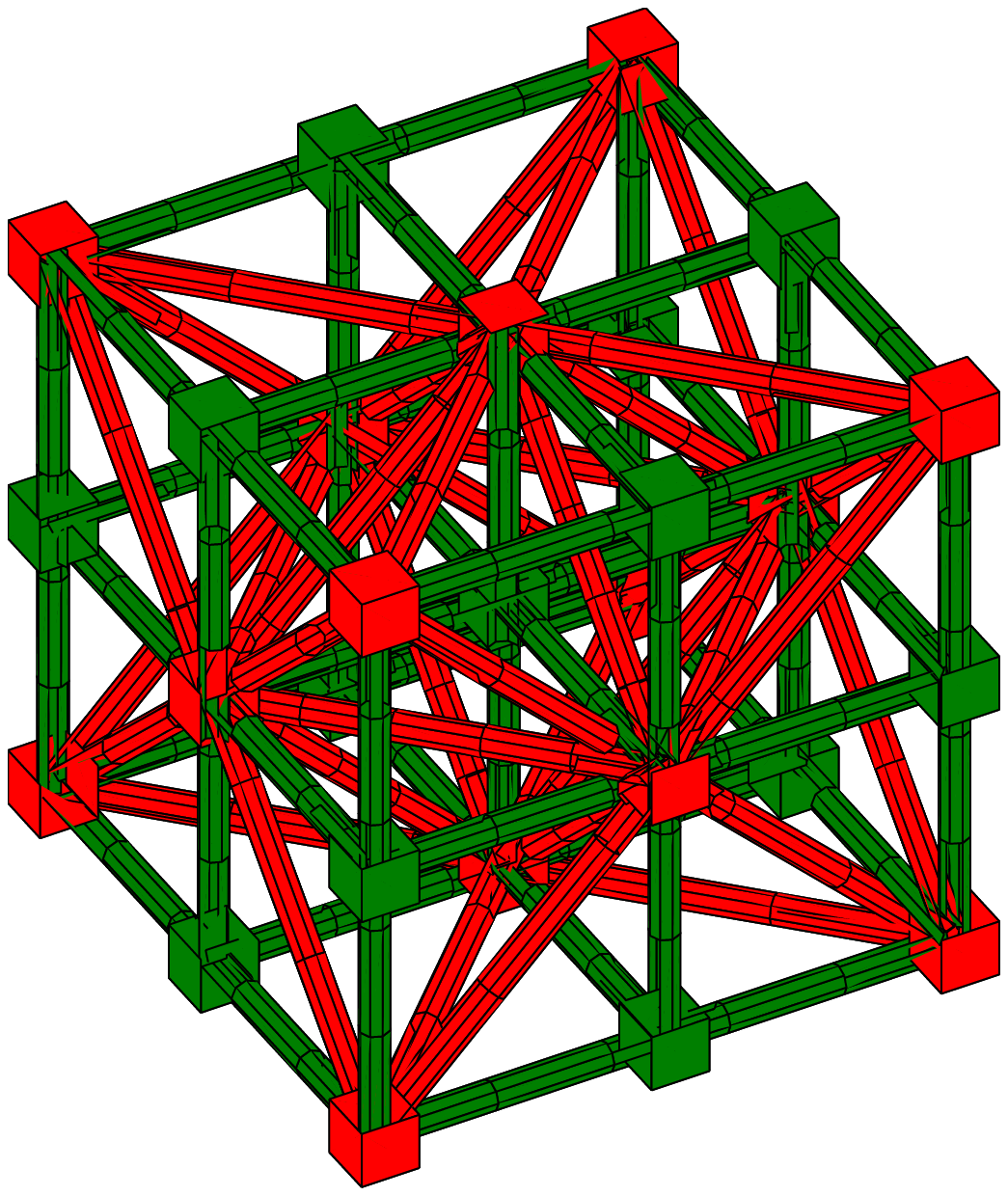}
\caption{the green and red grids superimposed showing the nodes
of the red grid at the corners and faces of the cube, and their
relationship to their six neighbouring nodes on the finer green grid.}
		\label{Fggr}
	\end{figure}
	The nodes of the red grid are at the corners of the cube and the
	centre of each of the faces as seen in Figure~\ref{Fggr}.  They
	each have six neighbours on the green grid so the natural
	restriction averaging of the residuals onto the red grid is
	\begin{eqnarray}
		r_{i,j,k}^{\ell-1}&=&\frac{1}{12}\left( 6r_{i,j,k}^\ell
		+r_{i+1,j,k}^\ell +r_{i-1,j,k}^\ell +r_{i,j+1,k}^\ell
+r_{i,j-1,k}^\ell
		+\right.\nonumber\\&&\left.\quad{}
		+r_{i,j,k+1}^\ell +r_{i,j,k-1}^\ell \right)\,,
		\label{Erred}
	\end{eqnarray}
	for $(i,j,k)$ corresponding to the (red) corners and faces of the
	coarse (blue) grid.  When the fine green grid is $n\times n\times
	n$ so that there are $N=n^3$ unknowns on the fine green grid, this
	average takes $8\flops$ for each of the approximately $N/2$ red
	nodes.  This operation count totals $4N\flops$.

	Note that throughout this discussion of restriction from the green
	to blue grids via the red and magenta, we index variables using
	subscripts appropriate to the fine green grid.  This also holds for
	the subsequent discussion of the prolongation from blue to green grids.

	\item
	\begin{figure}[tbp]
		\centering
		\includegraphics[width=0.9\textwidth]{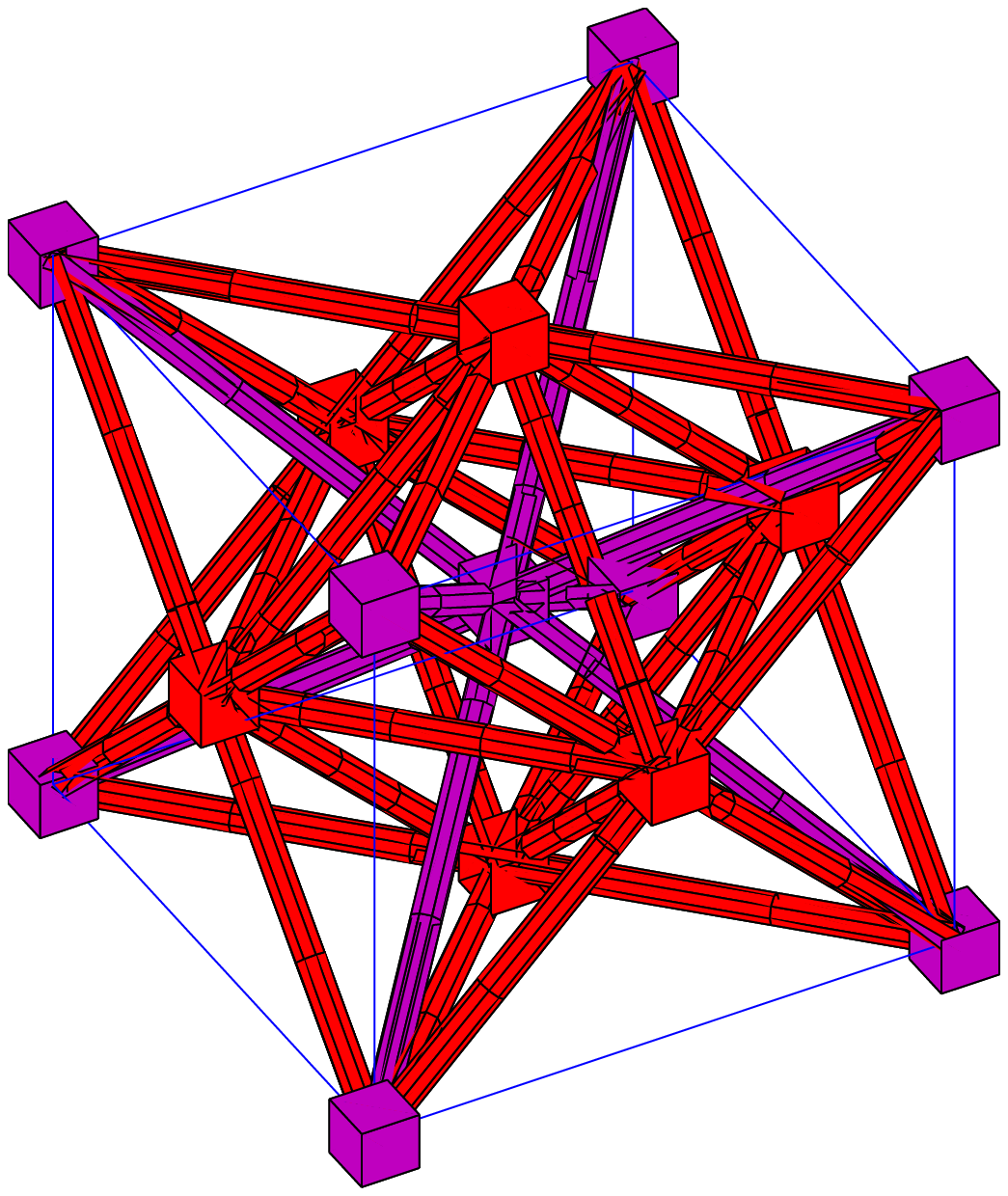}
		\caption{the red and magenta grids superimposed showing the
		nodes of the magenta grid at the corners and the centre of the
		(blue) cube.}
		\label{Frmag}
	\end{figure}
	The nodes of the next coarser grid, magenta, are at the corners and
	centres of the cube as seen in Figure~\ref{Frmag}.  Observe that the
	centre nodes of the magenta grid are not also nodes of the finer red
	grid; this causes some complications in the treatment of the two
	types of nodes.  The magenta nodes at the corners are connected to
twelve
	neighbours on the red grid so the natural average of the residuals
	is
	\begin{eqnarray}
		r_{i,j,k}^{\ell-1}&=&\frac{1}{24}\left( 12r_{i,j,k}^\ell
		+r_{i+1,j+1,k}^\ell +r_{i+1,j-1,k}^\ell
		+r_{i-1,j-1,k}^\ell +r_{i-1,j+1,k}^\ell
		+\right.\nonumber\\&&\left.\quad{}
		+r_{i+1,j,k+1}^\ell +r_{i+1,j,k-1}^\ell
		+r_{i-1,j,k-1}^\ell +r_{i-1,j,k+1}^\ell
		+\right.\nonumber\\&&\left.\quad{}
		+r_{i,j+1,k+1}^\ell +r_{i,j+1,k-1}^\ell
		+r_{i,j-1,k-1}^\ell +r_{i,j-1,k+1}^\ell
		\right)\,,
		\label{Ermagc}
	\end{eqnarray}
	for $(i,j,k)$ corresponding to the magenta corner nodes.  This
	average takes $14\flops$ for each of $N/8$ nodes.  The magenta
	nodes at the centre of the coarse (blue) cube is not connected to
	red nodes by the red grid, see Figure~\ref{Frmag}.  However, it
	has six red nodes in close proximity, those at the centre
	of the faces, so the natural average is
	\begin{equation}
		r_{i,j,k}^{\ell-1}=\frac{1}{6}\left( r_{i+1,j,k}^\ell
+r_{i-1,j,k}^\ell
		+r_{i,j+1,k}^\ell +r_{i,j-1,k}^\ell +r_{i,j,k+1}^\ell
+r_{i,j,k-1}^\ell
		\right)\,,
		\label{Ermagm}
	\end{equation}
	for $(i,j,k)$ corresponding to the magenta centre nodes.  This
	averaging takes $6\flops$ for each of $N/8$ nodes.  The operation
	count for all of this restriction step from red to magenta is
	$2\frac{1}{2}N\flops$.

	\item
	\begin{figure}[tbp]
		\centering
		\includegraphics[width=0.9\textwidth]{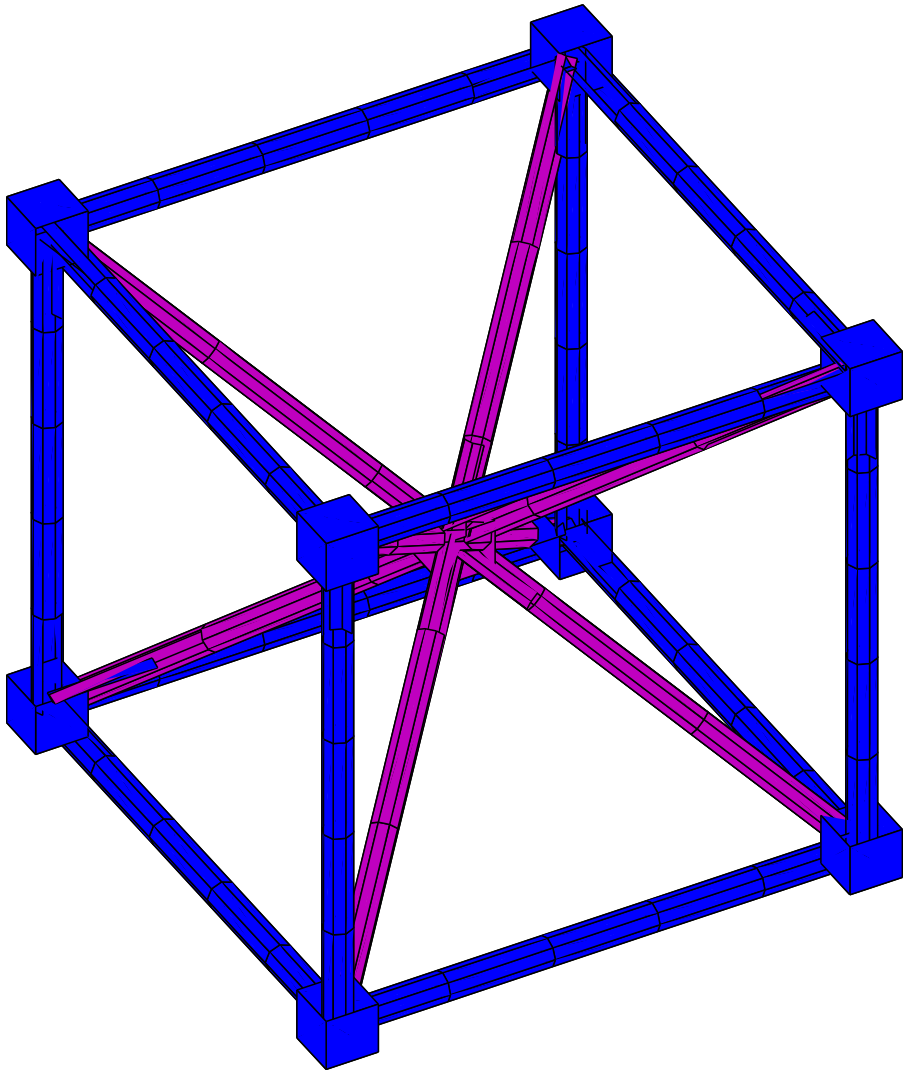}
		\caption{the magenta and blue grids superimposed showing the
		common nodes at the corners of the blue grid and the
connections to
		the magenta centre node.}
		\label{Frblu}
	\end{figure}
	The nodes of the coarse blue grid are at the corners of the shown
	cube, see Figure~\ref{Frblu}.  On the magenta grid they are
	connected to eight neighbours, one for each octant, so the natural
	average of residuals from the magenta to the blue grid is
	\begin{eqnarray}
		r_{i,j,k}^{\ell-1}&=&\frac{1}{16}\left( 8r_{i,j,k}^\ell
		+r_{i+1,j+1,k+1}^\ell +r_{i+1,j+1,k-1}^\ell
		+r_{i+1,j-1,k+1}^\ell
		+\right.\nonumber\\&&\left.\quad{}
		+r_{i+1,j-1,k-1}^\ell
		+r_{i-1,j+1,k+1}^\ell +r_{i-1,j+1,k-1}^\ell
		+\right.\nonumber\\&&\left.\quad{}
		+r_{i-1,j-1,k+1}^\ell +r_{i-1,j-1,k-1}^\ell
		 \right)\,,
		\label{Erblu}
	\end{eqnarray}
	for $(i,j,k)$ corresponding to the blue corner nodes.  This
	averaging takes $10\flops$ for each of $N/8$ blue nodes which thus
	totals $1\frac{1}{4}N\flops$.
\end{itemize}
These three restriction steps, to go up three levels of grids, thus
total approximately $7\frac{3}{4}N\flops$.  Hence, the entire
restriction process, averaging the residuals, from a finest grid of
spacing $h$ up $3L$ levels to the coarsest grid of spacing $H=2^Lh$
takes
\begin{equation}
	K_r\approx\frac{62}{7}N\left( 1-\frac{1}{8^L} \right)\flops \approx
	{\textstyle 8\frac{6}{7}}N\flops\,.
	\label{Ekrest3}
\end{equation}
The simplest standard one-step restriction direct from the fine green
grid to the blue grid takes approximately $3\frac{3}{4}N\flops$.
Over the whole hierarchy this totals $4\frac{2}{7}N\flops$ which is
roughly half that of the proposed method.  We anticipate that rapid
convergence of the V-cycle makes the increase worthwhile.

\subsection{The Jacobi prolongation steps}
\label{SSjp3}

As in 2D, with this rich structure of grids we have no need to
interpolate when prolongating from a coarse grid onto a finer grid; an
appropriate ``red-black'' Jacobi iteration of the residual
equation~(\ref{Evpois}) avoids interpolation.  Given an estimate of
corrections $v_{i,j,k}^\ell$ at some blue level grid we proceed to the
finer green grid via the following three prolongation steps.
\begin{itemize}
	\item Perform a magenta-blue Jacobi iteration on the nodes of the
	magenta grid shown in Figure~\ref{Frblu}.  See that each node on
	the magenta grid is connected to eight neighbours distributed
	symmetrically about it, each contributes to an estimate of the
	Laplacian at the node.  Thus, given initial approximations on the
	blue nodes from the coarser blue grid,
	\begin{eqnarray}
		v_{i,j,k}^\ell&=&\frac{1}{8}\left( -4p_mh^2r_{i,j,k}^\ell
		+v_{i+1,j+1,k+1}^{\ell-1} +v_{i+1,j+1,k-1}^{\ell-1}
		+v_{i+1,j-1,k+1}^{\ell-1}
		+\right.\nonumber\\&&\left.\quad{}
		+v_{i+1,j-1,k-1}^{\ell-1}
		+v_{i-1,j+1,k+1}^{\ell-1} +v_{i-1,j+1,k-1}^{\ell-1}
		+\right.\nonumber\\&&\left.\quad{}
		+v_{i-1,j-1,k+1}^{\ell-1} +v_{i-1,j-1,k-1}^{\ell-1}
		 \right)\,,
		\label{Emprolm}
	\end{eqnarray}
	for $(i,j,k)$ on the centre magenta nodes.  The following blue-Jacobi
	iteration uses these updated values in the similar formula
	\begin{eqnarray}
		v_{i,j,k}^\ell&=&\frac{1}{8}\left( -4p_mh^2r_{i,j,k}^\ell
		+v_{i+1,j+1,k+1}^{\ell} +v_{i+1,j+1,k-1}^{\ell}
		+v_{i+1,j-1,k+1}^{\ell}
		+\right.\nonumber\\&&\left.\quad{}
		+v_{i+1,j-1,k-1}^{\ell}
		+v_{i-1,j+1,k+1}^{\ell} +v_{i-1,j+1,k-1}^{\ell}
		+\right.\nonumber\\&&\left.\quad{}
		+v_{i-1,j-1,k+1}^{\ell} +v_{i-1,j-1,k-1}^{\ell}
		 \right)\,,
		\label{Emprolb}
	\end{eqnarray}
	for $(i,j,k)$ on the corner blue nodes.  In these formulae the
	over relaxation parameter $p_m$ has been introduced for later fine
	tuning; initially take $p_m=1$.  The operation count for this
	magenta-blue Jacobi iteration is $10\flops$ on each of $N/4$ nodes
	giving a total of $2\frac{1}{2}N\flops$.

	\item Perform a red-magenta Jacobi iteration on the nodes of the
	red grid shown in Figure~\ref{Frmag}.  However, because the centre
	node (magenta) is not on the red grid, two features follow: it is
	not updated in this prolongation step; and it introduces a little
	asymmetry into the weights used for values at the nodes.  The red
	nodes in the middle of each face are surrounded by four magenta
	nodes at the corners and two magenta nodes at the centres of the
	cube.  However, the nodes at the centres are closer and so have twice
	the weight in the estimate of the Laplacian.  Hence, given initial
	approximations on the magenta nodes from the coarser grid,
	\begin{eqnarray}
		v_{i,j,k}^{\ell}&=&\frac{1}{8}\left( -2p_{r1}h^2r_{i,j,k}^\ell
		+2\left[v_{i,j,k+1}^{\ell-1}+v_{i,j,k-1}^{\ell-1}\right]
		+\right.\nonumber\\&&\left.\quad{}
		+v_{i+1,j+1,k}^{\ell-1} +v_{i+1,j-1,k}^{\ell-1}
		+v_{i-1,j-1,k}^{\ell-1} +v_{i-1,j+1,k}^{\ell-1}
		\right)\,,
		\label{Erprolr}
	\end{eqnarray}
	for $(i,j,k)$ corresponding to the red nodes on the centre of
	faces normal to the $z$-direction.  Similar formulae apply for red
	nodes on other faces, cyclically permute the role of the indices.
	The over relaxation parameters $p_{r1}$ and $p_{r2}$ are
	introduced for later fine tuning; initially take
	$p_{r1}=p_{r2}=1$.  The following magenta-Jacobi iteration uses
	these updated values.  Each magenta corner node in
	Figure~\ref{Frmag} is connected to twelve red nodes and so is
	updated according to
	\begin{eqnarray}
		v_{i,j,k}^{\ell}&=&\frac{1}{12}\left(
		 -4p_{r2}h^2r_{i,j,k}^\ell
		+\right.\nonumber\\&&\left.\quad{}
		+v_{i+1,j+1,k}^\ell +v_{i+1,j-1,k}^\ell
		+v_{i-1,j-1,k}^\ell +v_{i-1,j+1,k}^\ell
		+\right.\nonumber\\&&\left.\quad{}
		+v_{i+1,j,k+1}^\ell +v_{i+1,j,k-1}^\ell
		+v_{i-1,j,k-1}^\ell +v_{i-1,j,k+1}^\ell
		+\right.\nonumber\\&&\left.\quad{}
		+v_{i,j+1,k+1}^\ell +v_{i,j+1,k-1}^\ell
		+v_{i,j-1,k-1}^\ell +v_{i,j-1,k+1}^\ell
		\right)\,,
		\label{Erprolm}
	\end{eqnarray}
	for all $(i,j,k)$ corresponding to corner magenta nodes.
	The operation count for this red-magenta Jacobi iteration
	is $9\flops$ on each of $3N/8$ nodes and $14\flops$ on each
	of $N/8$ nodes.  These total $5\frac{1}{8}N\flops$.

	\item Perform a green-red Jacobi iteration on the nodes of
	the fine green grid shown in Figure~\ref{Fggr}.  The green
	grid is a standard rectangular grid so the Jacobi
	iteration is also standard.  Given initial approximations on
	the red nodes from the coarser red grid,
	\begin{eqnarray}
		v_{i,j,k}^{\ell}&=&\frac{1}{6}\left( -p_gh^2r_{i,j,k}^\ell
		+v_{i+1,j,k}^{\ell-1} +v_{i-1,j,k}^{\ell-1}
		+v_{i,j+1,k}^{\ell-1} +v_{i,j-1,k}^{\ell-1}
		+\right.\nonumber\\&&\left.\quad{}
		+v_{i,j,k+1}^{\ell-1} +v_{i,j,k-1}^{\ell-1} \right)\,,
		\label{Egprolg}
	\end{eqnarray}
	for $(i,j,k)$ corresponding to the green nodes (edges and centre
	of the cube).  The over relaxation parameter $p_g$, initially
	$p_g=1$, is introduced for later fine tuning.  The red-Jacobi
	iteration uses these updated values in the similar formula
	\begin{eqnarray}
		v_{i,j,k}^{\ell}&=&\frac{1}{6}\left( -p_gh^2r_{i,j,k}^\ell
		+v_{i+1,j,k}^{\ell} +v_{i-1,j,k}^{\ell}
		+v_{i,j+1,k}^{\ell} +v_{i,j-1,k}^{\ell}
		+\right.\nonumber\\&&\left.\quad{}
		+v_{i,j,k+1}^{\ell} +v_{i,j,k-1}^{\ell} \right)\,,
		\label{Egprolr}
	\end{eqnarray}
	for the red nodes in Figure~\ref{Fggr}.  This prolongation
	step is a standard Jacobi iteration and takes $8\flops$ on
	each of $N$ nodes for a total of $8N\flops$.
\end{itemize}
These three prolongation steps together thus total
$15\frac{5}{8}N\flops$.  To prolongate over $\ell=3L$ levels
from the coarsest grid of spacing $H=2^Lh$ to the finest grid thus takes
\begin{equation}
	K_p\approx\frac{125}{7}N\left(1-\frac{1}{8^L}\right)\flops \approx
	{\textstyle 17\frac{6}{7}}N\flops\,.
	\label{Ekprol3}
\end{equation}
The simplest trilinear interpolation direct from the blue grid to the
green grid would take approximately $3\frac{1}{4}N\flops$, to be
followed by $8N\flops$ for a Jacobi iteration on the fine green grid.
Over the whole hierarchy this standard prolongation takes
approximately $12\frac{6}{7}N\flops$.  This total is smaller, but the
proposed diagonal grid achieves virtually three Jacobi
iterations instead of one and so is more effective.

\subsection{The V-cycle converges well}

Numerical experiments show that, as in 2D, although the operation
count of the proposed algorithm is a little higher, the speed of
convergence is much better.  Both algorithms use $9N\flops$ to compute
second-order accurate residuals on the finest grid.  Thus the proposed
method takes approximately $35\frac{5}{7}N\flops$ for one V-cycle,
some 37\% more than the $26\frac{1}{7}N\flops$ of the simplest
standard algorithm.   It achieves a mean factor of convergence
$\bar\rho\approx0.140$.  This rapid rate of convergence easily
compensates for the small increase in computations taking half the
number of flops per decimal digit accuracy determined.

\begin{table}[tbp]
	\centering
	\caption{comparison of cost, in flops, and performance for
	unoptimised
	algorithms for solving Poisson's equation in three spatial
	dimensions on a $17^3$ grid.  The column headed ``per iter'' shows
the number of
	flops per iteration, whereas column showing ``per dig'' is
	$\flops/\log_{10}\bar\rho_0$ and indicates the number of flops needed
	to compute each decimal digit of accuracy.}
	\begin{tabular}{|l|r|lr|}
		\hline
		algorithm & per iter & $\bar\rho_0$ & per dig   \\
		\hline
		diagonal, $\Ord{h^2}$ & $35.7N$ & 0.140 & $42N$  \\
		usual, $\Ord{h^2}$ & $26.1N$ & 0.477 & $81N$ \\
		\hline
		diagonal, $\Ord{h^4}$ & $48.7N$ & 0.659 & $269N$ \\
		usual, $\Ord{h^4}$ & $39.1N$ & 0.651 & $210N$ \\
		\hline
	\end{tabular}
	\label{tbl:3d}
\end{table}

As in 2D, fourth-order accurate solvers may be obtained simply by
using the above second-order accurate V-cycle on the fourth-order
accurate residuals evaluated on the finest grid.  A compact
fourth-order accurate scheme for the residuals is the 19~point
formula
\begin{eqnarray}
r_{i,j,k}&=&\frac{1}{12}\left( 6f_{i,j,k}
           +f_{i+1,j,k} +f_{i,j+1,k} +f_{i-1,j,k} +f_{i,j-1,k}
           +f_{i,j,k+1}
          +\right.\nonumber\\&&\quad\left.{}
          +f_{i,j,k-1}
           \right)
       -\frac{1}{6h^2}\left[ -24 u_{i,j,k}
         +2\left(u_{i,j-1,k} +u_{i,j+1,k} +u_{i-1,j,k}
          +\right.\right.\nonumber\\&&\quad\left.\left.{}
         +u_{i+1,j,k}
         +u_{i,j,k+1} +u_{i,j,k-1} \right)
	     +u_{i+1,j+1,k} +u_{i-1,j+1,k}
          +\right.\nonumber\\&&\quad\left.{}
	     +u_{i-1,j-1,k} +u_{i+1,j-1,k}
	     +u_{i,j+1,k+1} +u_{i,j+1,k-1} +u_{i,j-1,k-1}
          +\right.\nonumber\\&&\quad\left.{}
	     +u_{i,j-1,k+1}
	     +u_{i+1,j,k+1} +u_{i-1,j,k+1} +u_{i-1,j,k-1} +u_{i+1,j,k-1}
	      \right]\,.
	\label{Efos4}
\end{eqnarray}
Then using the V-cycle described above to determine corrections $v$ to
the field $u$ leads to an increase in the operation count of
$13N\flops$ solely from the extra computation in finding the finest
residuals.  Numerical experiments show that the multigrid iteration
still converge, albeit slower, with $\bar\rho\approx 0.659$.
Table~\ref{tbl:3d} shows that the rate of convergence on the diagonal
hierarchy of grids is little different than that for the simplest
usual multigrid algorithm.  As in 2D, high accuracy, 4th order
solutions to Poisson's equation are best found by employing a first
stage that finds 2nd order accurate solutions which are then refined
in a second stage.

\subsection{Optimise parameters of the V-cycle}
\label{SSopt3}

As in 2D, the multigrid algorithms are improved by introducing some
relaxation in the Jacobi iterations.  The four parameters $p_m$,
$p_{r1}$, $p_{r2}$ and $p_g$ were introduced in the Jacobi iterations
(\ref{Emprolm}--\ref{Egprolr}) to do this, values bigger than 1
correspond to some over relaxation.

\begin{table}[tbp]
	\centering
	\caption{comparison of cost, in flops, and performance for
	optimised algorithms for solving Poisson's equation in three
	spatial dimensions on a $17^3$ grid varying over relaxation
	parameters to determine the best rate of convergence.  The column
	headed ``per iter'' shows the number of flops per iteration,
	whereas column showing ``per dig'' is $\flops/\log_{10}\bar\rho$
	and indicates the number of flops needed to compute each decimal
	digit of accuracy.}
	\begin{tabular}{|l|r|lllllr|}
		\hline
		algorithm & per iter &  $p_{m}$ & $p_{r1}$
		& $p_{r2}$ & $p_{g}$ & $\bar\rho$ &
		per dig  \\
		\hline
		diag, $\Ord{h^2}$ & $35.7N$ & 1.11 & 1.42 & 1.08 &
		0.99 & 0.043 & $26N$  \\
		usual, $\Ord{h^2}$ & $26.1N$ &  &  &  & 1.30 & 0.31 &
		$51N$  \\
		\hline
		diag, $\Ord{h^4}$ & $48.7N$ & 0.91 & 0.80 & 0.70 &
		1.77 & 0.39 & $119N$  \\
		usual, $\Ord{h^4}$ & $39.1N$ &  &  &  & 1.70 & 0.41 &
		$101N$  \\
		\hline
	\end{tabular}
	\label{tbl:3dopt}
\end{table}

The search for the optimum parameter set used the Nelder-Mead simplex
method encoded in the procedure \textsc{fmins} in \matlab{}.  Searches
were started from optimum parameters found for coarser grids.  As
tabulated in Table~\ref{tbl:3dopt} the optimum parameters on a $17^3$
grid\footnote{Systematic searches on a finer grid were infeasible
within one days computer time due to the large number of unknowns:
approximately 30,000 components occur in the eigenvectors on a $33^3$
grid.} were $p_m=1.11$, $p_{r1}=1.42$, $p_{r2}=1.08$ and $p_g=0.99$
and achieve an astonishingly fast rate of convergence of
$\bar\rho\approx 0.043$.  This ensures convergence to a specified
precision at half the cost of the similarly optimised, simple
conventional multigrid algorithm.

For the fourth-order accurate residuals an optimised diagonal
multigrid performs similarly to the optimised conventional multigrid
with a rate of convergence of $\bar\rho\approx 0.39$.  Again fourth
order accuracy is best obtained after an initial stage in which second
order accuracy is used.

\section{Conclusion}

The use of a hierarchy of grids at angles to each other can halve the
cost of solving Poisson's equation to second order accuracy in grid
spacing.  Each iteration of the optimised \emph{simplest} multigrid
algorithm decreases errors by a factor of at least 20.  This is true
in both two and three dimensional problems.  Further research is
needed to investigate the effective of extra Jacobi iterations at each
level of the diagonal grid.

When compared with the amazingly rapid convergence obtained for the
second order scheme, the rate of convergence when using the fourth
order residuals is relatively pedestrian.  This suggests that a
multigrid V-cycle specifically tailored on these diagonal grids for
the fourth order accurate problem may improve convergence markedly.

There is more scope for W-cycles to be effective using these diagonal
grids because there are many more levels in the multigrid hierarchy.
An exploration of this aspect of the algorithm is also left for
further research.

\paragraph{Acknowledgement:} This research has been
supported by a grant from the Australian Research Council.

\bibliographystyle{plain}\bibliography{ajr,bib,new}

\begin{thebibliography}{1}

\bibitem{Dendy97}
J.E. Dendy.
\newblock Semicoarsening multgrid for systems.
\newblock {\em Elect. Trans. Num. Anal.}, 6:97--105, 1997.

\bibitem{Gupta97}
M.M. Gupta, J.~Kouatchou, and J.~Zhang.
\newblock Comparison of second- and fourth-order discretizations for multigrid
  {P}oisson solvers.
\newblock {\em J.~Comput. Phys.}, 132:226--232, 1997.

\bibitem{Ibraheem96}
S.O. Ibraheem and A.O. Demuren.
\newblock On bi-grid local mode analysis of solution techniques for {3-D Euler
  and Navier-Stokes} equations.
\newblock {\em J. Comput. Physics}, 125:354---377, 1996.

\bibitem{Iyengar90}
R.K. Iyengar and A.~Goyal.
\newblock A note on multigrid for the three-dimensional {P}oisson equation in
  cylindrical coordinates.
\newblock {\em J.~Comput. \& Appl. Math.}, 33:163--169, 1990.

\bibitem{Mavriplis97b}
D.J. Mavriplis.
\newblock Directional coarsening and smoothing for anisotropic {Navier-Stokes}
  problems.
\newblock {\em Elect. Trans. Num. Anal.}, 6:182--197, 1997.

\bibitem{Mavriplis97}
D.J. Mavriplis.
\newblock Unstructured grid techniques.
\newblock {\em Annu.\ Rev.\ Fluid Mech.}, 29:473--514, 1997.

\bibitem{Moulton98}
J.D. Moulton, J.E. Dendy, and J.M. Hyman.
\newblock The black box multigrid numerical homogenization algorithm.
\newblock {\em J. Comput. Physics}, 142:80---108, 1998.

\bibitem{Sheffer98}
S.G. Sheffer, L.~Martinelli, and A.~Jameson.
\newblock An efficient multigrid algorithm for compressible reactive flows.
\newblock {\em J. Comput. Physics}, 144:484---516, 1998.

\bibitem{Zhang98}
Jun Zhang.
\newblock Fast and high accuracy multigrid solution of the three-dimensional
  {P}oisson equation.
\newblock {\em J.~Comput.\ Phys.}, 143:449--461, 1998.

\end{thebibliography}

\end{document}